\documentclass[%
 aip,
reprint,
 amsmath,amssymb,
 aps,
showkeys
]{revtex4-1}

\usepackage{graphicx}
\usepackage{dcolumn}
\usepackage{bm}
\usepackage{appendix}
\usepackage{xr}
\usepackage[normalem]{ulem}
\usepackage[english]{babel}

\usepackage{todonotes}

\begin{document}

\preprint{CHA}

\title{Analyzing Collective Motion with Machine Learning and Topology}

\author{Dhananjay Bhaskar}
\affiliation{%
 Center for Biomedical Engineering, Brown  University,\\ Providence, RI, USA, 02912
}%

\author{Angelika Manhart}
\affiliation{%
 Department of Mathematics, University College London, \\ London, UK, WC1E 6BT
}%

\author{Jesse Milzman}
\affiliation{%
 Department of Mathematics, University of Maryland, \\ College Park, MD, USA, 20742
}%

\author{John T. Nardini}
\affiliation{%
 The Statistical and Applied Mathematical Sciences Institute (SAMSI),\\  Durham, NC, USA, 27709
}

\author{Kathleen M. Storey}
\affiliation{%
 Department of Mathematics, University of Michigan,\\  Ann Arbor, MI, USA, 48105
}%

\author{Chad M. Topaz}
\affiliation{%
 Department of Mathematics and Statistics, Williams College,\\ Williamstown, MA, USA, 01267
}%

\author{Lori Ziegelmeier}
\email{lziegel1@macalester.edu }
\affiliation{%
 Department of Mathematics, Statistics, and Computer Science, Macalester College, \\ Saint Paul, MN, USA, 55105
}%

\date{\today}

\begin{abstract}

We use topological data analysis and machine learning to study a seminal model of collective motion in biology [D'Orsogna \emph{et al.}, Phys. Rev. Lett. 96 (2006)]. This model describes agents interacting nonlinearly via attractive-repulsive social forces and gives rise to collective behaviors such as flocking and milling. To classify the emergent collective motion in a large library of numerical simulations and to recover model parameters from the simulation data, we apply machine learning techniques to two different types of input. First, we input time series of order parameters traditionally used in studies of collective motion. Second, we input measures based in topology that summarize the time-varying persistent homology of simulation data over multiple scales. This topological approach does not require prior knowledge of the expected patterns. For both unsupervised and supervised machine learning methods, the topological approach outperforms the one that is based on traditional order parameters.
\end{abstract}

\maketitle

{\bf A fundamental goal in the study of complex, nonlinear systems is to understand the link between local rules and collective behaviors. For instance, what influence does coupling between oscillators have on the ability of a network to synchronize? How does policing affect hotspots of crime in urban areas? Why do the movement decisions that fish make lead to schooling structures? We examine such links by bringing together tools from applied topology and machine learning to study a seminal model of collective motion that replicates behavior observed in biological swarming, flocking, and milling. Studies of collective motion often focus on the so-called forward problem: given a particular mathematical model, what dynamics are observed for different parameters input to the model? In contrast, we study an inverse problem: given observed data, what model parameters could have produced it? Also, given observed data, what paradigmatic dynamics are being exhibited? To answer these questions, we use machine learning techniques and find that they achieve higher accuracy when applied to topological summaries of the data -- called crockers -- as compared to more traditional summaries of the data that are commonly used in biology and physics to characterize collective motion. Ours is the first study to use crockers for nonlinear dynamics classification and parameter recovery.}

\section{\label{sec:intro}Introduction}

Fundamental to many nonlinear systems is the link between local rules and global behaviors. For instance, what influence does coupling between oscillators have on the ability of a network to synchronize? How does policing affect hotspots of crime in urban areas? Why do the movement decisions that fish make lead to schooling structures? In this paper, we consider the relationship between local rules and global behavior in service of two tasks of interest in nonlinear science, namely classification of dynamics and parameter recovery. We conduct this study in the context of one of the  aforementioned applications: collective motion in biology.

Animal groups such as flocks, herds, schools, and swarms are ubiquitous in nature \cite{CamDenFra2001,Sum2010}, inspire numerous mathematical models\cite{giardina2008,vicsek2012,degond2018,degond2018age}, and motivate biomimetic approaches to engineering and computer science.\cite{VinBogBog2006,Bhu2009} In animal groups, two levels of behavior come into play. First, at the individual level, organisms make decisions about how to move through space. It is well-established in the biological literature that social interactions between organisms play a key role. Second, at the group level, collective motion may occur, where animals coordinate and produce emergent patterns. Mathematical models can help link these two levels of behavior. This type of linkage is fundamental to innumerable nonlinear systems displaying collective behavior.


While there is a vast literature on mathematical models of collective motion, we focus on the influential model of D'Orsogna \emph{et al.}\cite{DorChuBer2006} This model describes agents whose movement is determined by self-propulsion, drag, and social attraction-repulsion, forces frequently used in collective motion studies. \cite{LevRapCoh2001,CouKraJam2002} The model produces paradigmatic behaviors such as rotating rings, vortices, disorganized swarms, and  traveling groups; these mimic fish schools, swarms of midges, bird flocks, and more. \cite{ParEde1999,Hep1997,NiOue2015}

Studies of collective motion models often focus on a \emph{forward problem}: given a particular model, what dynamics are observed for different parameters? In our present work, however, we study an \emph{inverse problem}: given observed data, what parameters could have produced it? Also, given observed data, what paradigmatic dynamics are being exhibited?\cite{lukeman2010, manhart2018} There are numerous strategies to infer parameters from data with simple differential equation models, including Bayesian inference \cite{stuart_inverse_2010} and frequentist approaches \cite{banks_parameter_1983}. Our primary contribution here is to solve an inverse problem with an agent-based model using tools from topological data analysis and machine learning. Though we focus on collective motion in the D'Orsogna model, our protocol is applicable in many other settings.

Machine learning and mathematical modeling have traditionally been viewed as separate ways of understanding data. On one hand, machine learning can extract predictions of complex relations within large data sets. On the other hand, modeling can be used to hypothesize how mechanisms lead to an observed behavior. There is a growing understanding, however, that modeling and machine learning can be used in synergy. \cite{BakPenJay2018} For example, Lu \emph{et al}. use nonparamteric estimators to learn the rules governing the observed output of agent-based models.\cite{lu_nonparametric_2019} The method is applied to fundamental interacting particle systems \cite{jones_determination_1924}, models of social influence \cite{couzin_effective_2005}, predator-swarm systems \cite{chen_minimal_2014}, and phototactic bacteria \cite{ha_particle_2009}.

Our study hinges on the use of topological data analysis (TDA), a set of tools that ``help the data analyst summarize and visualize complex datasets.''\cite{Was2018} TDA has played a pivotal role in studies of breast cancer, spinal cord injury, contagion, and other biological systems. \cite{NicLevCar2011,NiePaqLiu2015,TayKliHar2015} A primary tool in the TDA toolbox is \emph{persistent homology}. \emph{Homology} has to do with calculating certain topological characteristics, while \emph{persistence} refers to examining which of these are maintained across multiple scales in the data. In its fundamental form, persistent homology provides a framework for describing the topology of a static data set. However, ideas from persistent homology have been extended to time-varying data. One approach is the \emph{Contour Realization Of Computed $k$-dimensional hole Evolution in the Rips complex}, known more simply as a \emph{crocker}. \cite{TopZieHal2015} A crocker shows contours of quantities called \emph{Betti numbers} as a function of time and of persistence scale, providing a topological summary of time-varying point clouds of data. Recent work has shown how crockers can be used for exploratory data analysis of collective motion and for judging the fitness of potential mathematical models of experimental data. \cite{TopZieHal2015,UlmZieTop2019}

Our present work brings together mathematical modeling, machine learning, and TDA to study collective motion in the D'Orsogna model. More specifically, we construct crockers from numerical simulations of the model and use them as inputs to machine learning clustering and classification algorithms in order to identify different paradigmatic patterns and the model parameters that produce them. We compare this approach to a more traditional one, which uses order parameters commonly used in studies of collective motion. While the traditional order parameters are typically chosen using prior knowledge of the system, the TDA tools can be used with no prior knowledge and are problem-independent. In our methodological study, we use data simulated from models so that we can compare inferred parameters to those actually used to generate the data. Once trained, the machine learning algorithms can be used to infer model parameters from experimental data.

The rest of this paper is organized as follows. Section~\ref{sec:methods} describes our methods, namely, numerical simulation of the D'Orsogna model, computation of traditional order parameters and crockers, and machine learning techniques. Section~\ref{sec:results} presents our results, including our primary finding: the topological approach outperforms the traditional one. We conclude in Section~\ref{sec:discussion}.

\section{\label{sec:methods}Methods}

\subsection{\label{sec:dorsogna}D'Orsogna Model and Numerical Simulations}

The D'Orsogna model\cite{DorChuBer2006,ChuDorMar2007} describes the motion of $N$ interacting agents of mass $m$. Each agent is characterized by its position and velocity $\boldsymbol{x}_i, \boldsymbol{v}_i \in \mathbb{R}^d$, $i = 1,\ldots,N$. We focus on the two-dimensional case, $d=2$, throughout this study, though the case for $d=3$ has been considered in Ref.~\onlinecite{chuang_swarming_2016}.  Position and velocity obey the coupled, nonlinear ordinary differential equations:
\begin{subequations}
\label{eq:dorsogna}
\begin{align}
\dot{\boldsymbol{x}}_{i} & =\boldsymbol{v}_{i},\label{eq:dorsogna1}\\
m\dot{\boldsymbol{v}}_{i} & =(\alpha-\beta|\text{\textbf{v}}_{i}|^{2})\text{\textbf{v}}_{i}-\vec{\nabla}_{i}U(\text{\ensuremath{\boldsymbol{x}}}_{i}),\label{eq:dorsogna2}\\
U(\boldsymbol{x}_{i}) & =\sum_{j\ne i}^{N}\left[C_{r}e^{-|\boldsymbol{x}_{i}-\boldsymbol{x}_{j}|/\ell_{r}}-C_{a}e^{-|\boldsymbol{x}_{i}-\boldsymbol{x}_{j}|/\ell_{a}}\right]\label{eq:dorsogna3},
\end{align}
\end{subequations}
where $\vec{\nabla}_{i}$ is the gradient with respect to $\boldsymbol{x}_i$. Eq.~\eqref{eq:dorsogna1}  states that the time derivative of position is velocity. Eq.~\eqref{eq:dorsogna2} is Newton's law, with the right hand side describing three forces acting on each agent: self-propulsion with strength $\alpha$, nonlinear drag with strength $\beta$, and social interactions. These social interactions are attractive-repulsive, as specified by the Morse-type potential in \eqref{eq:dorsogna3}. The first term inside the brackets describes social repulsion of overall strength $C_r > 0$. Repulsion decays exponentially in space, with characteristic length scale of decay $\ell_r > 0$. The exponential decay reflects that organisms' ability to sense each other through sight, sound, or smell decays over distance. Restated, an organism will be more heavily influenced by near individuals than far ones. The second term is signed oppositely, and hence describes social attraction. The summation in \eqref{eq:dorsogna3} means that a given organism interacts with all other organisms, albeit with influence decaying exponentially in space.

Biological modeling studies typically assume that repulsion is strong but operates over a short length scale, while attraction is weaker but operates over a longer scale. For \eqref{eq:dorsogna}, this would mean $C_r > C_a$ and $\ell_r < \ell_a$. In this case, the Morse potential $U$ has a well-defined minimum, representing the distance at which attraction and repulsion balance for two organisms in isolation. However, this distance is not the separation observed between individuals in a group because of the nonlinear all-to-all coupling. Regardless, as in Ref.~\onlinecite{DorChuBer2006}, we do not enforce the restriction $C_r > C_a$ and $\ell_r < \ell_a$. 

For our study, we set $N=200$. After nondimensionalizing \eqref{eq:dorsogna} we are left with four dimensionless parameters: $\alpha$, $\beta$, $C=C_r/C_a$ and $l=l_r/l_a$. We set $\alpha=1.5$, and $\beta=0.5$, corresponding to a base case of parameters from Ref.~\onlinecite{DorChuBer2006}, and explore the remaining parameter space by varying the ratios $C = C_r/C_a$ and $\ell = \ell_r/\ell_a$. Both $C$ and $\ell$ will take on values in $\{0.1,0.5,0.9,2.0,3.0\}$, resulting in 25 different possible parameter combinations. For each combination, we perform 100 simulations using M{\sc ATLAB}'s \texttt{ode45} function. Because the D'Orsogna model is typically independent of initial conditions, we draw $\mathbf{x}_{i}(0)$ and $\mathbf{v}_{i}(0)$ each from a uniform distribution on $[-1,1]^d$.  We simulate until $t=100$, allowing the swarm to attain a dynamic equilibrium state. From the computed trajectories, we sample the positions and velocities every 0.05 time units. Thus, the final output is $\{(\boldsymbol{x}_i(t_j),\boldsymbol{v}_i(t_j))\}_{i=1,..,N}^{j=1,...,M}$ for $t_j= (j-1)\Delta t, M=2001, \Delta t =  0.05$ for each of our 2500 simulations.

\begin{figure*}
    \centering
    \includegraphics[width=.9\textwidth]{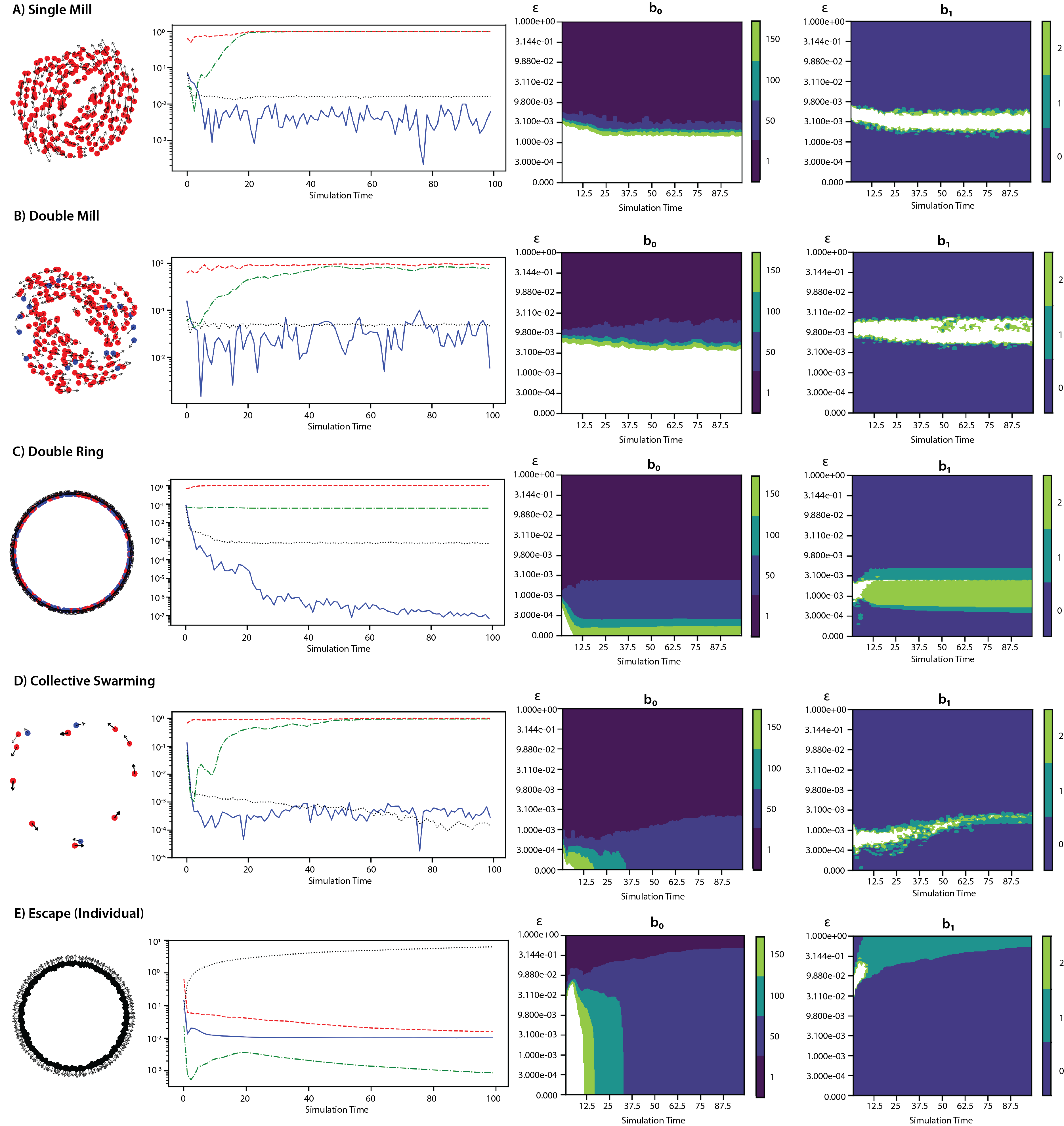}
    \caption{\label{fig:phenotypes} Collective motion phenotypes  generated by the D'Orsogna model \eqref{eq:dorsogna} with $N=200$ agents, $\alpha=1.5$, and $\beta=0.5$. First column (not to scale): Snapshots of agents' positions (dots) and velocities (arrows) at $t=100$ with clockwise motion in blue and counter-clockwise motion in red (where applicable). Second column: Order parameter time series, namely  polarization $P(t)$ (blue), angular momentum $M_{ang}(t)$ (green), absolute angular momentum $M_{abs}(t)$ (red), and average distance to nearest neighbor $D_{NN}(t)$ (black). Third and fourth columns: time-delay crockers showing, respectively, Betti numbers $b_0$ and $b_1$ as a function of time $t$ and persistence parameter $\epsilon$ (log scale). For ease of depiction, we represent any value of $b_0>150$ or $b_1>2$ as white. We include only a small subset of contours for visual clarity. White regions correspond to larger values of Betti numbers, not shown. (A) Single mill, $C=0.5$, $\ell=0.1$. (B) Double mill, $C=0.9$, $\ell=0.5$. (C) Double ring, $C=0.1$, $\ell=0.1$. (D) Collective swarming, $C=0.1$, $\ell=0.5$. (E) Escape (individual), $C=2.0$, $\ell=0.9$.}
\end{figure*}

\begin{figure}[ht]
    \centering
    \includegraphics[width=0.49\textwidth]{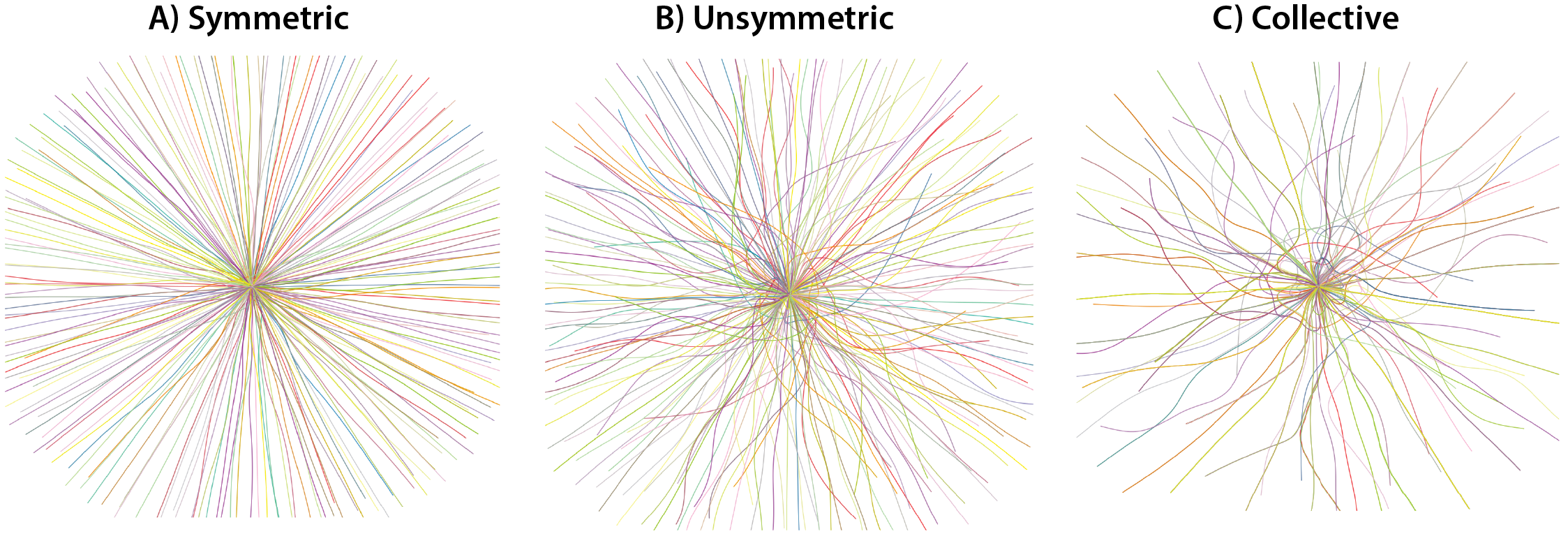}
    \caption{Three sub-classes of the escape phenotype observed in the D'Orsogna model \eqref{eq:dorsogna}. $N$, $\alpha$, and $\beta$ are as in Fig.~\ref{fig:phenotypes}. Each agent's trajectory is shown in a randomly chosen color for times $t \in [0, 40]$.
(A) $(C, \ell) = (2, 0.9)$. Particles escape to infinity individually, and the overall pattern is radially symmetric.
(B) $(C, \ell) = (2, 2)$. Agents escape individually, but form gaps in angular distribution.
(C) $(C, \ell) = (2, 3)$. Agents superpose in small groups and move outwards with large gaps in angular distribution. \label{fig:escape}}
\end{figure}

These simulations produce paradigmatic collective motion including single mills, double mills, double rings (referred to simply as rings in Ref.~\onlinecite{DorChuBer2006}), collective swarms, and group escape, which we split into three distinct classes. \cite{DorChuBer2006} For the remainder of this paper, we refer to paradigmatic collective behaviors as \emph{phenotypes}. The first column of Fig.~\ref{fig:phenotypes} shows a representative snapshot of each major phenotype, and Fig.~\ref{fig:escape} shows the three distinct escape types. Table~\ref{tab:phenotype_params} lists the values of $(C,\ell)$ that produce each phenotype in our library of simulations. We describe these phenotypes in more detail at the end of Section~\ref{subsec:tda}.

\begin{table}[h]
\caption{\label{tab:phenotype_params}Collective motion phenotypes and corresponding social interaction potential parameters in our library of simulations of the D'Orsogna model \eqref{eq:dorsogna}. Here, $C=C_r/C_a$, $\ell = \ell_r/\ell_a$, and $N=200$. We fix the remaining parameters $\alpha=1.5$ and $\beta=0.5$ as in Ref.~\onlinecite{DorChuBer2006}.}
\begin{ruledtabular}
\begin{tabular}{ll}
\textbf{Phenotype} & \textbf{Parameters} $\mathbf{(C,\ell)}$   \\\hline\hline
 Single mill & (0.5,0.1),(0.9,0.1),(2,0.1),(2,0.5), \\
 & (3,0.1)\\
 Double mill & (0.9,0.5)  \\
 Double ring & (0.1,0.1),(0.5,0.5),(0.9,0.9)  \\
 Collective swarm & (0.1,0.5),(0.1,0.9),(0.1,2),(0.1,3), \\
 & (0.5,0.9),(0.5,2),(0.5,3),(0.9,2),(0.9,3)\\
 Escape (symmetric) & (2,0.9),(3,0.9)\\
 Escape (unsymmetric) & (2,2),(3,2),(3,3)\\
 Escape (collective) & (2,3),(3,0.5)\\
\end{tabular}
\end{ruledtabular}
\end{table}

\subsection{Order Parameters}\label{subsec:op}

Investigators often use \emph{order parameters} to summarize the output of collective motion experiments or simulations. \cite{CouKraJam2002,DorChuBer2006,TopZieHal2015} Summaries are necessary because it is impractical or impossible to manually inspect large amounts of raw output data. The order parameters are intended to suggest if and when certain types of group behavior emerge in a population, for instance, when many individuals move in the same direction or rotate with the same orientation. Typical order parameters used for \eqref{eq:dorsogna} include polarization, angular momentum, absolute angular momentum, and the mean distance to the nearest neighbor\cite{TopZieHal2015}. We use these four in our present study, plotted for each phenotype in the second column of Fig.~\ref{fig:phenotypes}.  

Group polarization $P(t)$ measures the degree of alignment between agents and is given by
\begin{align}
P(t) = \left| \frac{\sum_{i=1}^N \textbf{v}_i(t)}{\sum_{i=1}^N |\textbf{v}_i(t)|} \right| \in [0,1],
\end{align}
with $P=1$ signifying that all agents have the same direction of motion. All phenotypes in Fig.~\ref{fig:phenotypes} exhibit low $P(t)$, suggesting no translational flocking.
Angular momentum $M_{ang}(t)$ can detect rotational motion, and is given by
\begin{align}
M_{ang}(t) = \left| \frac{\sum_{i=1}^N \textbf{r}_i(t) \times\textbf{v}_i(t)}{\sum_{i=1}^N |\textbf{r}_i(t)| |\textbf{v}_i(t)|} \right|  \in [0,1],
\end{align}
where $\textbf{r}_i(t) = \textbf{x}_i(t) - \textbf{x}_{cm}(t)$, and $\textbf{x}_{cm}(t)$ refers to the center of mass of the agents. A group with $M_{ang}=1$ would have individuals sharing perfectly rotational motion. Following Ref.~\onlinecite{ChuDorMar2007}, we also consider the absolute angular momentum $M_{abs}(t)$, given by
\begin{align}
M_{abs}(t) = \left| \frac{\sum_{i=1}^N | \textbf{r}_i(t) \times\textbf{v}_i(t) |}{\sum_{i=1}^N |\textbf{r}_i(t)| |\textbf{v}_i(t)|} \right|  \in [0,1].
\end{align}
 Discrepancies between angular momentum and absolute angular momentum can distinguish a single mill from a double mill, in which counter-rotating agents cancel out each others' angular momentum. For example, we observe in Fig.~\ref{fig:phenotypes}  that $M_{\text{abs}}(t)\approx M_{\text{ang}}(t)$ for the single mill, whereas $M_{\text{abs}}(t) > M_{\text{ang}}(t)$ for the double mill.  Finally, we consider the mean distance to nearest neighbor, $D_{NN}(t)$, which may distinguish group escape behavior from other phenotypes. \cite{Huepe2008} $D_{NN}(t)$ is given by
\begin{align*}
D_{NN}(t) = \frac{1}{N}\displaystyle\sum_{i=1}^N \min_{j\neq i} |\textbf{x}_i(t) - \textbf{x}_j(t)| \in \mathbb{R}_{\geq 0}.
\end{align*}
 $D_{NN}(t)$ becomes very large for the escape phenotype over time in Fig.~\ref{fig:phenotypes} as the particles act repulsively.

In calculating time series of these order parameters, we downsample by a factor of 23 (chosen since it is a divisor of 2001, the original number of simulation frames), resulting in $M=87$ time points. While some information is lost with downsampling, it makes subsequent computations faster while maintaining a high level of classification accuracy with the downsampling rate we have chosen.

\subsection{Persistent Homology and Crockers}\label{subsec:tda}

While order parameters can be useful summaries of collective motion data, they are typically designed in a problem-specific manner and with some knowledge of the expected dynamics. We compare an order parameter approach to one that measures the topology of the data. This approach is arguably more agnostic and less application-specific. We now review relevant ideas from topology. To make this review broadly accessible, we keep it conceptual. For some technical details, see, \emph{e.g.}, Refs.~\onlinecite{Hat2002,TopZieHal2015,OttPorTil2017}.

We begin our explanation of persistent homology and crockers by focusing on agents' positions during one time step of a simulation (or experiment). This data constitutes a \emph{point cloud}, made up of $N$ points in $\mathbb{R}^d$. To study the topology of a point cloud, we transform it into an object called a \emph{simplicial complex}. While there are many ways to construct a simplicial complex, we use the \emph{Vietoris-Rips} (VR) complex, a common choice in TDA because it is efficient to compute.

To build a VR complex, we select a distance $\epsilon$ and draw a ball of diameter $\epsilon$ around each point. If two balls intersect, we connect them with an edge. If three balls all pairwise intersect, we connect all three edges and fill in the resulting triangle. If four balls all pairwise intersect, we connect all six edges, fill in each of the four triangles bounded by those edges, and fill in the solid tetrahedron bounded by the four triangles. Points, edges, filled triangles, and solid tetrahedra are called \emph{0-, 1-, 2-}, and \emph{3-simplices}, and more generally, $k$-simplices for any $k+1$ points with $\epsilon$-balls that pairwise intersect.

With a simplicial complex built from our point cloud, we now measure its topology by calculating its \emph{Betti numbers}. Betti numbers $b_k$ are topological invariants, meaning that they are unchanged under continuous deformations of the object such as stretching, compressing, warping, and bending. Thus, they measure something fundamental about the shape of the object. More specifically, Betti numbers enumerate the number of distinct holes in the complex that have a $k$-dimensional boundary, that is, a hole surrounded by $k$-simplices. For instance, $b_0$ counts the number of connected components in the simplicial complex. Similarly, $b_1$ counts the number of topological loops that bound a 2-D void. Betti number $b_2$ counts the number of trapped 3-D volumes, and so on as dimension increases. Algebraic topology tells us how to encode the calculation of Betti numbers as a linear algebra problem; see standard topology texts or Ref.~\onlinecite{Top2015} for a tutorial. The calculation is a \emph{homology} computation because $b_k$ is the rank of an algebraic object called a \emph{homology group}.

In the discussion above, no value of $\epsilon$ was specified. \emph{Persistent homology} constructs simplicial complexes and calculates Betti numbers for a range of $\epsilon$ values. There exist powerful software packages that automate this process. \cite{OttPorTil2017} We use the \texttt{Ripser} package. \cite{TraSauBar2018} The outputs of these computations are \emph{birth} and \emph{death} values of $\epsilon$, that is, the values of $\epsilon$ for which the various features enumerated by $b_k$ appear and disappear. The word \emph{persistence} refers to the ranges of $\epsilon$ over which features persist. For example, features that persist over large ranges of $\epsilon$ might be interpreted as signals rather than topological noise. There exist many ways to organize the birth and death information, with the most common being objects called barcodes and persistence diagrams. \cite{OttPorTil2017} 
Additionally, for a given value of $k$, one could construct a vector in which each entry gives the value of $b_k$ for a specific value of $\epsilon$ (say, on a grid). This information, $b_k(\epsilon)$, is a \emph{Betti curve}.

Small perturbations in data produce small perturbations in persistence diagrams; that is to say, persistence diagrams are stable to noise near the data.\cite{Cohen-Steiner2007} However, persistence computations are not stable with respect to outliers. In practice, a codensity measure can be used to filter outliers. Further work in this area develops notions of distance that are robust to large quantities of empirical noise and outliers.\cite{Chazal2011,Fasy2014, Chazal2018} A clustering approach could also be used to limit the effects of this type of noise.\cite{Diky2019}

Thus far, we have discussed topological analysis of static point clouds. If we allow our agents' positions to evolve dynamically in time, $t$, then we can construct a Betti curve for each frame of the simulation or experiment, and concatenate these into a matrix. We let time $t$ vary along columns and $\epsilon$ vary along rows.  Each entry specifies the Betti number $b_k$ for a specific pair $(t,\epsilon)$. The matrix is a topological signature of the time-varying data of a simulation and once vectorized, can serve as input to machine learning algorithms. Equivalently, for visualiation, one could take this matrix and construct a contour plot of $b_k(t,\epsilon)$. Such a plot is the \emph{Contour Realization Of Computed $k$-dimensional hole Evolution in the Rips complex}, or \emph{crocker}, defined in Ref. \onlinecite{TopZieHal2015}.

\begin{figure}[ht]
    \centering
    \includegraphics[width=0.48\textwidth]{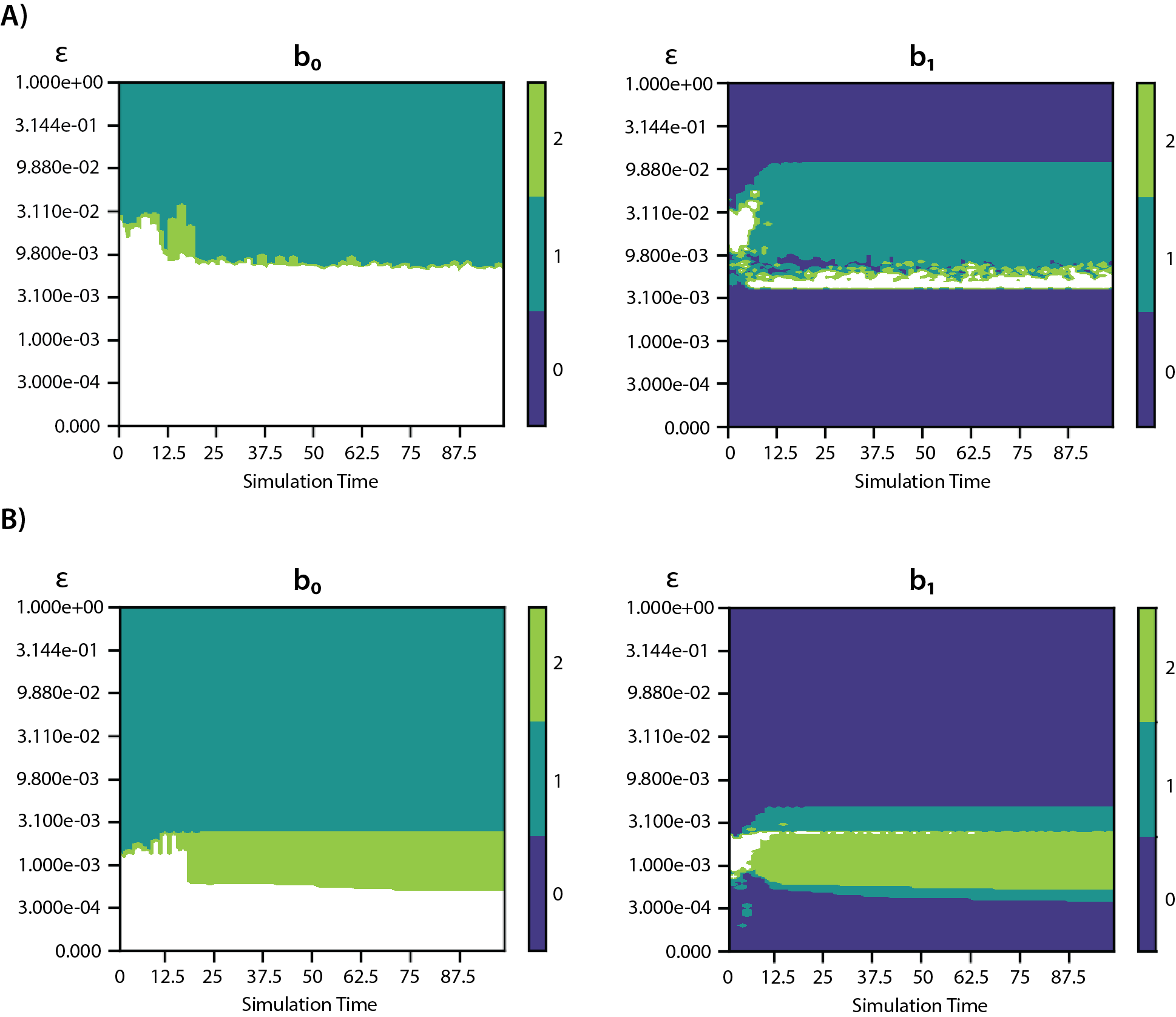}
    \caption{Crockers for $b_0$ (left) and $b_1$ (right) computed on a simulation of the D'Orsogna model \eqref{eq:dorsogna} with $(C, \ell) = (0.1,0.1)$. The simulation produces a double ring phenotype consisting of groups of agents moving clockwise and counterclockwise around a circle. For ease of depiction, we represent any contours greater than 2 as white. (A) Crockers obtained from position and velocity data. We normalize position and velocity by using their respective maximum and minimum magnitudes across all simulations. This approach identifies only a single topological loop. (B) Crockers computed from 4-D data that incorporates time-delayed position. This approach avoids a disparity between the scales of position and velocity, and correctly detects two loops (and two connected components).} \label{fig:crocker_plots}
\end{figure}

As mentioned in Section~\ref{sec:dorsogna}, our data consist of numerical solutions of \eqref{eq:dorsogna}. While the D'Orsogna model tracks agents' positions and velocities, we restrict ourselves to using position data in our topological analysis. This approach has two advantages. First, it renders our techniques applicable to experimental data, where position is the most easily observed quantity. Second, it circumvents a potential scaling disparity between numerical values of position and velocity when performing TDA on the data, as exhibited in Figure ~\ref{fig:crocker_plots}. In contrast, three out of four order parameters described above require knowledge of the velocity.

To regain some of the information lost by excluding velocity, we incorporate time-delayed position information into some of our analyses. We demonstrate this time-delay approach for a double ring phenotype in Fig.~\ref{fig:crocker_plots}, finding $b_1=2$ for a range of $\epsilon$ values, as we would expect based on Ref.~\onlinecite{TopZieHal2015}. Our time-delayed point cloud consists of points in $\mathbb{R}^4$ of the form $(\boldsymbol{x}_i(t_j), \boldsymbol{x}_i(t_j-5\Delta t))$, where $j$ ranges from 23 to 2001 with a spacing of 23 to result in 86 time samples.  Fig.~\ref{fig:4D_visual} and accompanying text in Appendix~\ref{app:4Dvis} describe the 4-D data for single mills, double mills, and double rings. We will also consider position-only crockers (computed on a point cloud consisting of $(\boldsymbol{x_i(t_j)})$, where the sampling of $j$ is the same as that discussed for the order parameters in Section ~\ref{subsec:op}).

Still, challenges remain with the normalization of our topological data. With escape phenotypes, inter-agent distances can approach infinity, whereas they remain bounded for other phenotypes. To circumvent this challenge, when performing our topological analyses, we take any agent whose distance from the origin crosses the threshold $\|\boldsymbol{x}_i\|_\infty = 10$ and edit the simulation data to hold it fixed at this position.

Even after enforcing this bound, phenotypes occur on a range of scales. While agent coordinate positions are capped at 10 for group escape, they are as small as $10^{-3}$ for collective swarms. We normalize position data across all simulations with a global normalization constant to ensure $-1\le \| \boldsymbol{x}_i(t)\|_\infty \le1$. With this scaling, the smallest normalized phenotypes have typical distances of $10^{-4}$. Thus, we compute persistent homology with $\epsilon$ varying logarithmically between $10^{-4}$ and 1 with 200 grid points such that $\epsilon_q= 10^{-4+q\Delta\epsilon}, \Delta\epsilon = 4/200, q = 1,...,200$. 

The third and fourth columns of Fig.~\ref{fig:phenotypes} show crockers ($b_0$ and $b_1$) for five example simulations. The crockers for single and double mills differ markedly. In (B), the $b_1$ crocker for the double mill contains small islands corresponding to two loops ($b_1=2$) within a large area of topological noise ($b_1>2$). In (A), the $b_1$ crocker for the single mill lacks this signature. In (C), a very strong signature of two loops ($b_1=2$) appears for the double ring simulation. Appendix~\ref{app:4Dvis} provides an explanation for the presence of two loops for double mills and rings. In (D), multiple agents form tight clumps in the collective swarm simulation, with each clump sufficiently dense that it appears as a single dot in the figure. The time-scale at which clumps form manifests as the disappearance of high-valued regions in the $b_0$ contour. On a macroscopic scale, we notice that each clump eventually travels with rotational motion, consistent with $b_1=1$ over a range of scales at later times. In (E), agents escape to infinity in a radially expanding circular arrangement. The strong signature of $b_1=1$ occurs at larger scales as time increases, consistent with an expanding circle.

\subsection{Unsupervised Learning}
\label{sec:unsupervised}

We use the $k$-medoids algorithm to cluster numerical simulations. Each simulation is characterized by a \emph{feature vector} constructed either from traditional order parameters or from topology. The order parameter feature vectors consist of time series $P(t)$, $M_{ang}(t)$, $M_{abs}(t)$, $D_{NN}(t)$, or the concatenation of all four. The topological feature vectors consist of vectorized crocker matrices for $b_0$, $b_1$, or the concatenation of the two, calculated from agent position (in some cases, augmented with time-delayed position as described in Section~\ref{subsec:tda}).

The $k$-medoids algorithm divides the ensemble of simulations into $k$ clusters, each of which is defined by one member of the ensemble that serves as the \emph{medoid}. The algorithm chooses medoids to minimize the sum of pairwise distances within each cluster, and each simulation is assigned to the cluster containing its closest medoid. We use the R software function \texttt{pam} to cluster our simulations into $k=25$ groups, since there are 25 distinct parameter choices ($C$,$\ell$). This is an \emph{unsupervised} approach, as it does not require labeled training data.

\subsection{Supervised Learning}
\label{sec:supervised} 

As an alternative approach, we use a multiclass linear support vector machine (SVM) to infer parameters \cite{cortes_support-vector_1995}. Our use of SVMs is \emph{supervised} because we train them on a subset of our simulations, each labeled with its true $(C,\ell)$ values. A linear SVM takes this training data and finds hyperplanes that maximally divide the simulations according to parameter values. To classify a simulation not included in the training set, one identifies the intra-hyperplane region in which it falls and reads off the appropriate label, \emph{i.e.}, the parameter values.

We use Matlab's \texttt{fitcecoc} function to build and train our SVMs using a one-versus-one approach and 5-fold cross validation. That is, for each round of cross-validation, we withhold 20\% of the data (20 simulations from each parameter combination) for the testing data set. We then train the linear SVM on the remaining data and compute the out-of-sample accuracy for simulations in the testing set.

Order parameter and topological feature vectors are not of the same dimension. The time series of each order parameter is 87-dimensional, and the time series of all four concatenated is $4 \times 87=348$-dimensional. On the other hand, each position crocker is $200\times 87=17400$-dimensional, and the concatenation of $b_0$ and $b_1$ is double that. To make a fair comparison between the order parameter and topological approaches, we use principal component analysis\cite{Jol1986} (PCA)  to reduce the dimensionality of our input feature vectors. In one case, we reduce crockers and the concatenated order paramers to 87 dimensions in order to compare them directly to the individual order parameter time series. In a second case, we reduce all feature vectors to three dimensions to investigate performance at low dimensionality. 

\section{\label{sec:results}Results}

\begin{figure}[h!]
    \centering
    \includegraphics[width=0.48\textwidth]{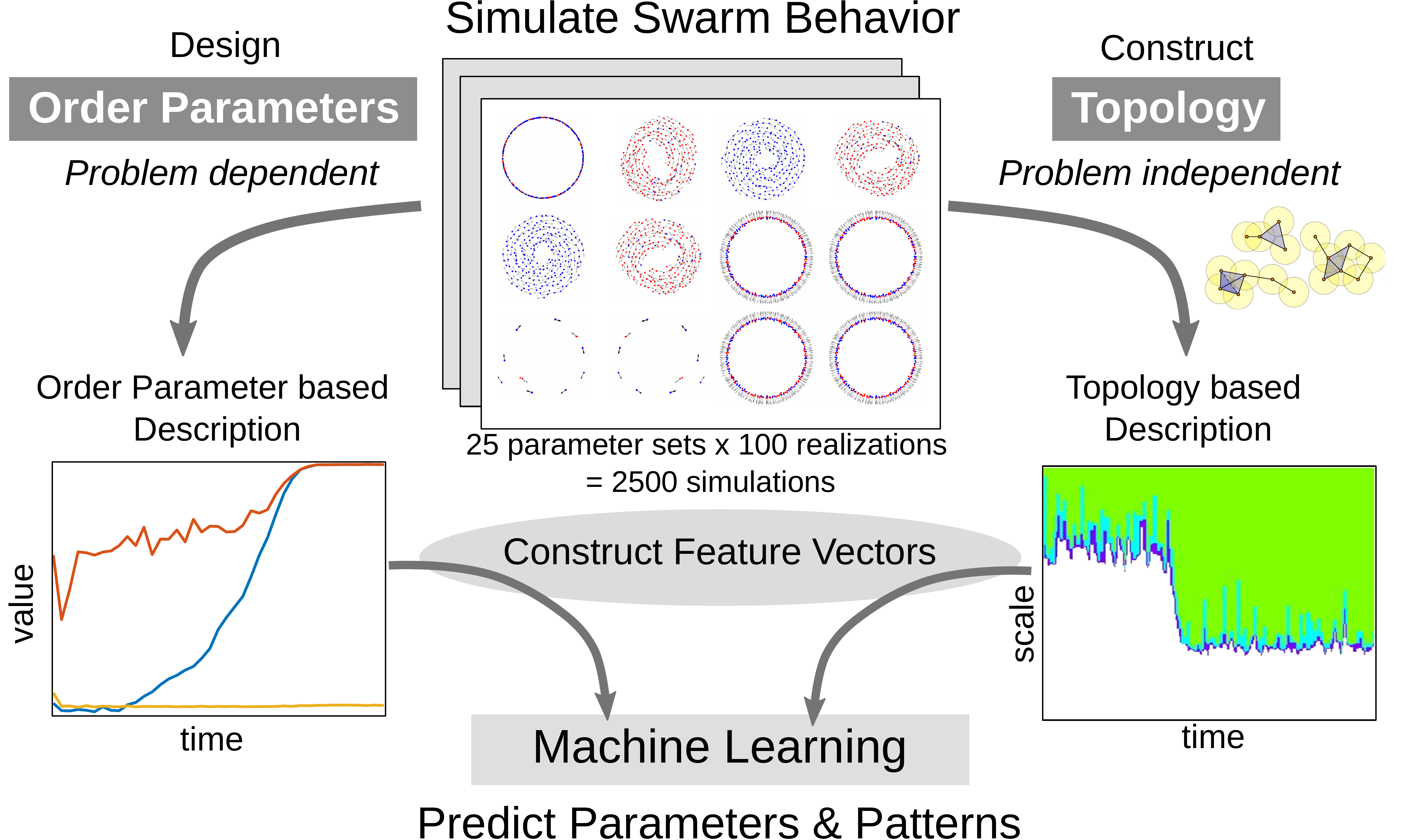}
    \caption{Our analysis pipeline. We summarize the dynamics of the D'Orsogna model \eqref{eq:dorsogna} using problem-specific order parameters (left) and a problem-independent description based on topology (right). We construct feature vectors from each summary and input them into machine learning algorithms to identify parameters and phenotype for each model simulation.}
    \label{fig:graphSummary}
\end{figure}

 Fig.~\ref{fig:graphSummary} recapitulates the entire analysis pipeline for our study. We summarize the dynamics of \eqref{eq:dorsogna} using two approaches. The order parameter-based approach uses problem-specific quantities designed to distinguish between observed dynamics. The topology-based approach does not require this \emph{a priori} knowledge and is instead based purely on the shape of the data. We then construct feature vectors from each type of summary and input them into machine learning algorithms to identify parameters and a phenotypic pattern for each model simulation. In the following subsections, we compare how accurately these different approaches can recover simulation parameters and identify dynamic phenotypes.

\subsection{Unsupervised Learning Results\label{sec:unsupervised}}

Table~\ref{tab:unsup_results} summarizes results for $k$-medoids clustering with $k=25$. Columns 1 and 2 specify the feature vector used, while columns 3 and 4 give the accuracies obtained for parameter recovery and phenotype identification. For the parameter recovery task, using the concatenation of all four order parameters yields 49.9\% accuracy. For the concatenation of $b_0$ and $b_1$ based only on position data, we obtain 76.6\% accuracy. Finally, for $b_0$ and $b_1$ based on time-delayed position, we have 71.3\% accuracy. In Appendix~\ref{app:MLresults}, Fig.~\ref{fig:unsup_confusion} displays confusion matrices. These reveal that regardless of feature vector type, the collective swarming and single mills are the most frequently misclassified phenotypes.

\begin{table}
     \caption{\label{tab:unsup_results}Unsupervised classification accuracy for parameter values using $k$-medoid clustering with various input feature vectors and $k=25$. The third column displays the accuracy when simulations are classified by parameter vector $(C,\ell)$, and the fourth column displays the accuracy when simulations are classified by phenotype.}
    \begin{tabular}{|c|c|c|c|}
    \hline
    \textbf{Summary} & \textbf{Feature}  & \textbf{Parameter}  & \textbf{Phenotype}  \\ \hline
    \hline
          & $P(t)$ & 20.6\% & 62.9\%\\
          & $M_{ang}(t)$ & 19.3\% & 60.7\%\\
        Order   &  $M_{abs}(t)$ & 47.3\% & 63.9\%\\
         Parameters & $D_{NN}(t)$ & 48\% & 80.1\%\\
          & All & 49.9\% & 97.4\% \\ \hline
          TDA & $b_0$ & 76.6\% & 83.6\% \\
           (position)& $b_1$& 71.9\% & 79.8\%\\
         & $b_0$ \& $b_1$ & 76.6\% & 99.8\%\\ \hline
         TDA & $b_0$ & 71.1\% & 82.6\%\\
          (time-delayed & $b_1$& 73.7\% & 84.4\%\\
          position) & $b_0$ \& $b_1$ & 71.3\% & 99.0\%\\ 
        \hline  
    \end{tabular}
\end{table}

For all feature vector types, phenotype identification is significantly more accurate than parameter recovery. The confusion matrices reveal that most cases of parameter recovery failure nonetheless place a simulation in its appropriate phenotypic regime. Using topological feature vectors based on position, we classify nearly every simulation correctly by phenotype, and this approach slightly outperforms concatenation of all order parameters. The slight improvement in classification accuracy when using TDA instead of order parameters may not be sufficient to justify the increased computational time for phenotype classification tasks. However, for parameter recovery tasks, TDA significantly improves classification accuracy when using k-medoids, and, as discussed shortly, when using a supervised classification method.

\subsection{Supervised Learning Results\label{sec:supervised}}

Table \ref{tab:supervised_results} summarizes supervised classification results for linear SVMs. For topological feature vectors based on position, $b_0$ does best with 97.0\% accuracy. Similarly, for delayed positions, $b_0$ also does best, with 99.6\% accuracy. Finally, for order parameter feature vectors, $D_{NN}(t)$ does best with 91.1\% accuracy, while concatenating all four order parameters yields 89.2\% accuracy.

For any feature vectors with dimension greater than 87, we also include results obtained after reducing the dimensionality to 87 via a principal component analysis, allowing for a more fair comparison. After dimension reduction, the time-delayed topological information for $b_1$ achieves the highest  classification accuracy at 99.9\%, followed by the position-only topological information for $b_0$ at 96.2\%. The classification accuracy for  all four concatenated order parameters drops to 69.6\%. Thus, a four-fold reduction in the dimension of the concatenated order parameters results in a 19.6\% loss of accuracy, whereas a 200-fold reduction for time-delayed topological information leads merely to a 0.8\% drop in accuracy for position-only information and a 0.1\% increase for time-delayed topological information. These results suggest that even with dimensionality reduction, the topological feature vectors still carry more discriminative information. 

To examine the limit of low-dimensional data, we calculate accuracies obtained after reducing all feature vectors to three dimensions. In this case, for topological feature vectors based on position data, the concatenation of $b_0$ and $b_1$ does best, with an accuracy of 93.1\%. For delayed position data, $b_0$ does best, achieving $87.7\%$ accuracy, and for order parameters, $D_{NN}(t)$ does best, yielding 81.5\% accuracy. Fig.~\ref{fig:PCA_depiction} in Appendix~\ref{app:MLresults} shows the three-dimensional representations of the $b_0$ and $b_1$ crockers. We observe a strong separation of the different phenotypes, which explains the high out-of-sample classification accuracy.

Also in Appendix~\ref{app:MLresults}, Fig.~\ref{fig:supervised_plots} visualizes the classification results. These results suggest that, similarly to the unsupervised case, collective swarms are the most difficult phenotype to classify. Still, overall, using topological data rather than order parameter data can significantly improve parameter recovery.

\begin{table}[]
    \centering
    \begin{tabular}{|c|c|c|c|}
    \hline
    \textbf{Summary} & \textbf{Feature}  & \textbf{Dimension}  & \textbf{Accuracy}  \\ \hline \hline
          & $P(t)$ & 87 & 57.7\%\\
          & $M_{ang}(t)$ & 87 & 34.4\%\\
          & $M_{abs}(t)$ & 87 & 68.0\%\\
          & $D_{NN}(t)$ & 87 & 91.1\%\\
          & All & 4$\times$87 & 89.2\% \\
  Order        & All (PCA) & 87 & 69.6\% \\
Parameters     & $P(t)$ (PCA) & 3 & 46.7\%\\
          & $M_{ang}(t)$ (PCA) & 3 & 30.0\%\\
          & $M_{abs}(t)$ (PCA) & 3 & 58.8\%\\
          & $D_{NN}(t)$ (PCA) & 3 & 81.5\%\\
          & All (PCA) & 3 & 68.6\% \\ 
          
          \hline
          & $b_0$ & 200$\times$87 & 97.0\% \\
          & $b_1$& 200$\times$87 & 93.7\%\\
          & $b_0$ and $b_1$ & 2$\times$200$\times$87 & 96.4\%\\
          TDA & $b_0$ (PCA) & 87& 96.2\%\\
          (position) & $b_1$ (PCA) & 87& 95.2\%\\
          & $b_0$ \& $b_1$ (PCA) & 87 & 96.2\%\\ 
          & $b_0$ (PCA) & 3 & 93.0\%\\
          & $b_1$ (PCA) & 3& 79.4\%\\
          & $b_0$ \& $b_1$ (PCA) & 3& 93.1\%\\ 
          \hline
          & $b_0$ & 200$\times$86& 99.6\%\\
          & $b_1$& 200$\times$86 & 99.3\%\\
          & $b_0$ \& $b_1$ & 2$\times$200$\times$86 & 99.1\%\\
          TDA & $b_0$ (PCA) & 87 & 99.7\%\\
          (time-delayed & $b_1$ (PCA) & 87 & 99.9\%\\
           position) & $b_0$ \& $b_1$ (PCA) & 87 & 99.7\%\\ 
          & $b_0$ (PCA) & 3 & 89.7\%\\
          & $b_1$ (PCA) & 3 & 82.8\%\\
          & $b_0$ \& $b_1$ (PCA) & 3 & 89.6\%\\ 
          \hline
          
     \hline
    \end{tabular}
    \caption{Supervised classification accuracy for parameter recovery using a trained linear SVM with various input feature vectors. In some cases, the dimensionality of feature vectors has been reduced using PCA.}
    \label{tab:supervised_results}
\end{table}

\section{\label{sec:discussion}Conclusions and Discussion}

We have combined mathematical modeling, topological data analysis, and machine learning to study nonlinear dynamics and parameter inference in the D'Orsogna model of collective motion. More specifically, we simulated \eqref{eq:dorsogna}, summarized the data using traditional and topological descriptors, and input these summaries into unsupervised and supervised machine learning algorithms in order to recover model parameters and classify pattern phenotypes.

Our machine learning classifiers achieved higher accuracy when using topological feature vectors (namely, crockers) than when using feature vectors based on traditional order parameters. Since the crocker feature vectors have higher dimensionality than the order parameter ones, we sought a fair comparison by reducing them via PCA. In this case, the crocker approach still achieved better classification accuracy using a supervised approach. In fact, $b_0$ crockers generated from time-delayed position data produced a nearly perfect classification. The time-delayed crockers encode information on particle velocity, which appears explicitly in Equation \eqref{eq:dorsogna} and may thus aid the algorithm. The addition of $b_1$ information serves primarily to increase the dimensionality of the feature space and results in reduced, though still quite high, classification accuracy.

One limitation of the topological approach is the computational cost required to produce crockers. While an order parameter is scalar-valued at each time, a crocker is vector-valued. However, recent software improvements, including the development of the \texttt{Ripser} package, have led to a significant reduction in cost.  

A major advantage of using topological data summaries is that they do not require prior knowledge about the patterns resulting from model simulation. Order parameters, on the other hand, are typically developed to capture specific features of previously observed model behavior. We found that for the D'Orsogna model, topological approaches to phenotype classification and parameter recovery achieved higher accuracy than order parameters even though they do not incorporate knowledge of the model or its dynamics.

In future work, we would like to apply this approach to data from biological experiments or field observations. There is a scarcity of publicly available data describing real biological aggregation dynamics, so for the present, we have demonstrated our method on simulation data. Furthermore, we would like to extend our work to more complex settings, \emph{e.g.}, to the D'Orsogna model posed in three dimensions, in which dynamical transitions occur between distinct phenotypic regimes \cite{chuang_swarming_2016}. Finally, it would be useful to augment the model with noise and assess its effect by using our topological methods.

\begin{acknowledgments}
This material is based upon work supported by the National Science Foundation under Grant Number DMS 1641020. We wish to acknowledge the American Mathematical Society's Mathematical Research Community which brought our collaboration together and supported our work. LZ is partially supported by NSF grant CDS\&E-MSS-1854703. DB is funded by the National Cancer Institute IMAT Program (R21CA212932). JN is partially supported by NSF
grant DMS-1638521 to SAMSI. CT is supported by NSF grant DMS-1813752.
\end{acknowledgments}

\appendix


\section{Time-Delayed Crockers}
\label{app:4Dvis}

\begin{figure*}
    \includegraphics[width=\textwidth]{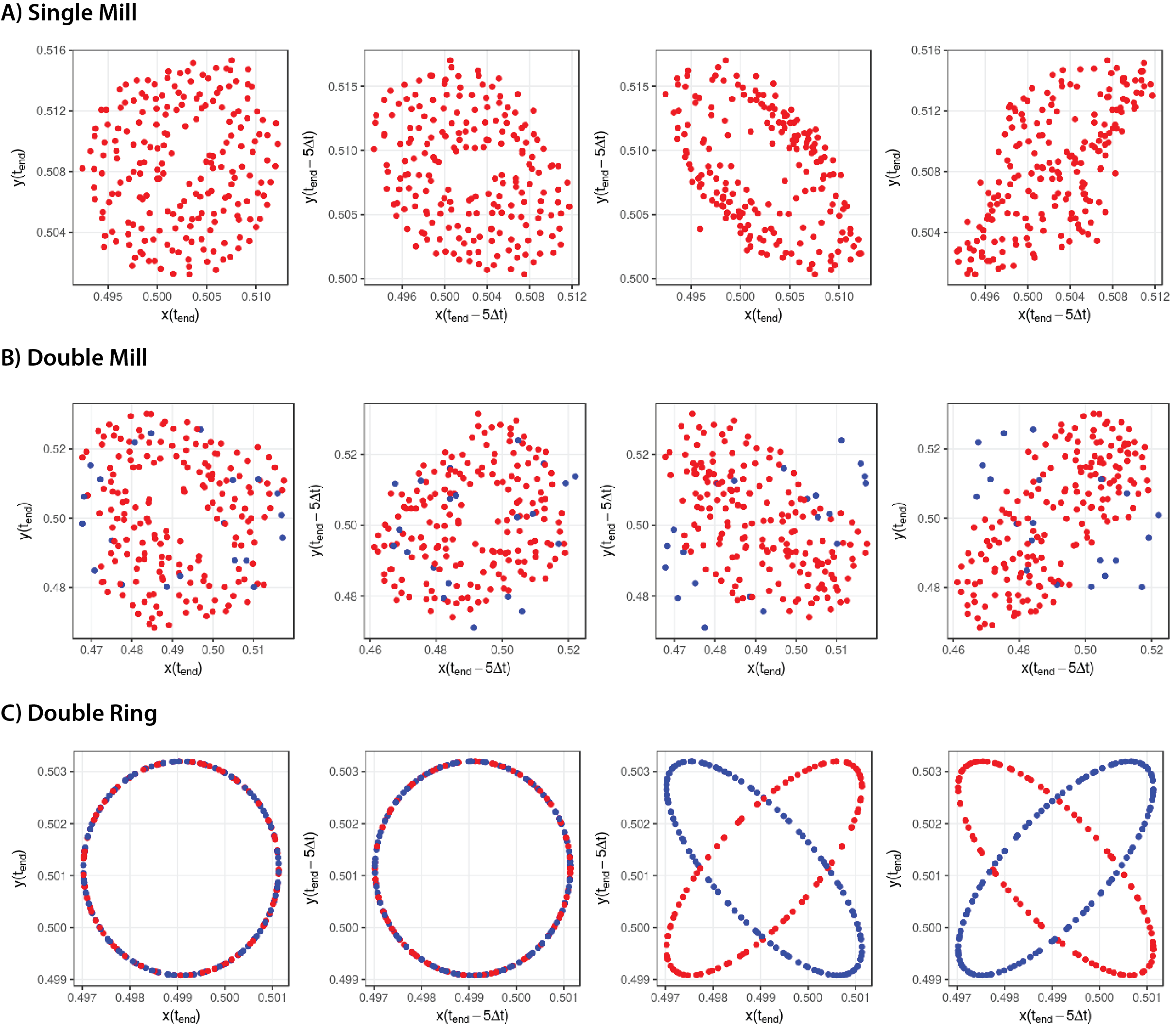}
    \caption{Snapshots ($t=100$) of data simulated from the D'Orsogna model \eqref{eq:dorsogna}, used to compute Betti numbers for the phenotypes in Fig.~\ref{fig:phenotypes}(A-C). Because the data is 4-D, as described in Section~\ref{subsec:tda}, we show various 2-D projections. Agents rotating clockwise around the group's center of mass are colored in blue, and agents moving counterclockwise are red. (A)~Single mill. (B)~Double mill. (C)~Double ring.}
    \label{fig:4D_visual}
\end{figure*}

In Fig.~\ref{fig:4D_visual}, we examine data used to compute crockers in order to understand differences in $b_1$ for the single mill, double mill, and double ring phenotypes. Since the data is 4-D, as described in Section~\ref{subsec:tda}, we show various 2-D projections. Panel (A) shows a single mill. We see a single loop-like structure that arises from the arrangement of particles in an annulus all traveling with the same orientation. This structure produces $b_1=1$ in Fig.~\ref{fig:phenotypes}(A). Panel (B) shows a double mill. Especially in the last two columns, we see a signal of two loops. The two loops correspond to the two counter-rotating swarms of the double mill. However, the signal is quite noisy. This noisy signal manifests as the transient islands of $b_1=2$ in Fig.~\ref{fig:phenotypes}(B).  Panel (C) shows a double ring. Agents occupy a well-defined circle, with some rotating clockwise and others counterclockwise. This phenotype gives rise to two loops in the 4-D space of our data (see last two columns) and produces a clear signal of $b_1=2$ in Fig.~\ref{fig:phenotypes}(C).

\section{Visualization of Machine Learning Results}
\label{app:MLresults}

\begin{figure}
    \centering
    \includegraphics[width=.5\textwidth]{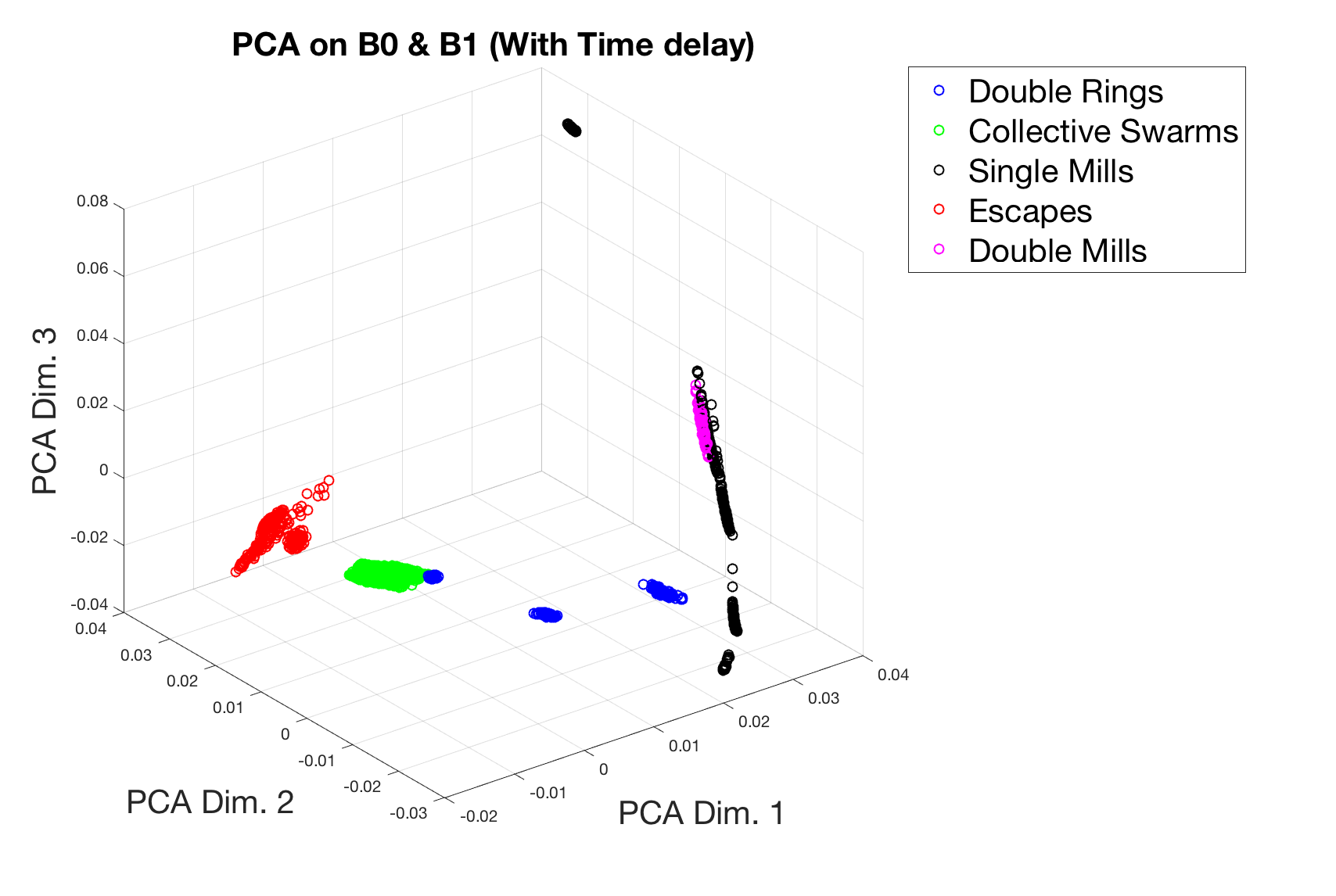}
    \caption{Depiction of the time-delayed $b_0$ and $b_1$ crockers after reduction to three dimensions with PCA.}
    \label{fig:PCA_depiction}
\end{figure}

\begin{figure}
    \centering
    \includegraphics[width=0.5\textwidth]{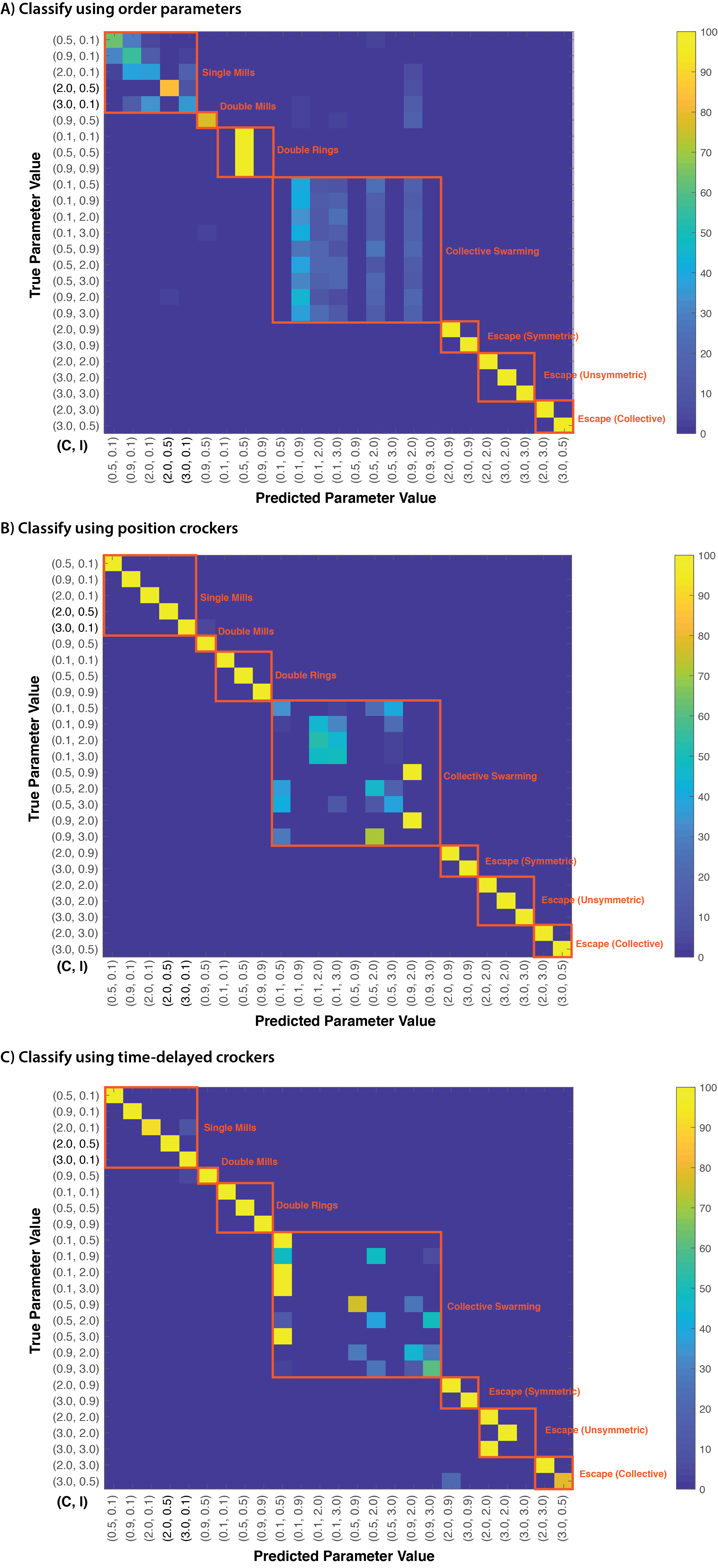}
    \caption{Confusion matrices for the unsupervised classification of Section~\ref{sec:unsupervised}. Using the $k$-medoids algorithm, we cluster simulations into $k=25$ groups. (A) Classification based on the concatenation of four order parameters. (B) Classification using the concatenation of $b_0$ and $b_1$ crockers computed from 2-D position data. (C) Like (B), but based on 4-D position data incorporating time delay. Rows correspond to actual parameter values ($C$,$\ell$) and columns correspond to those assigned by the classifier. There are 100 simulations for each true parameter bin. Color indicates the number of simulations for each true/predicted combination.}
    \label{fig:unsup_confusion}
\end{figure}

Fig.~\ref{fig:PCA_depiction} displays the representation of all 2500 time-delayed $b_0$ and $b_1$ crockers reduced to three-dimensional space using PCA. Here we see a strong separation of the single mill, double mill, collective swarm, escape, and double ring phenotypes. This separation explains why a linear SVM trained on this information has a fairly high out-of-sample classification accuracy, as shown in Table \ref{tab:supervised_results}.

Fig.~\ref{fig:unsup_confusion} displays confusion matrices arising from the unsupervised classification (k-medoids clustering) performed in Section~\ref{sec:unsupervised}. In panel (A), we cluster on the concatenation of all four order parameters. Most of the simulations with misidentified model parameters are those exhibiting single mill behavior or collective swarm behavior. However, for these incorrect cases, most were clustered with cases sharing the same phenotype. In panel (B), we use the concatenation of $b_0$ and $b_1$ crockers derived from 2-D position data. In this case, parameter recovery is more accurate for single mill simulations than in panel (A), but collective swarms remain challenging. Panel (C) is similar to (B), but is based on 4-D data incorporating position and time-delayed position. This approach yields results similar to those in (B), with a slight increase in misclassification among the three escape phenotypes. 

\begin{figure}
    \centering
    \includegraphics[width=0.4\textwidth]{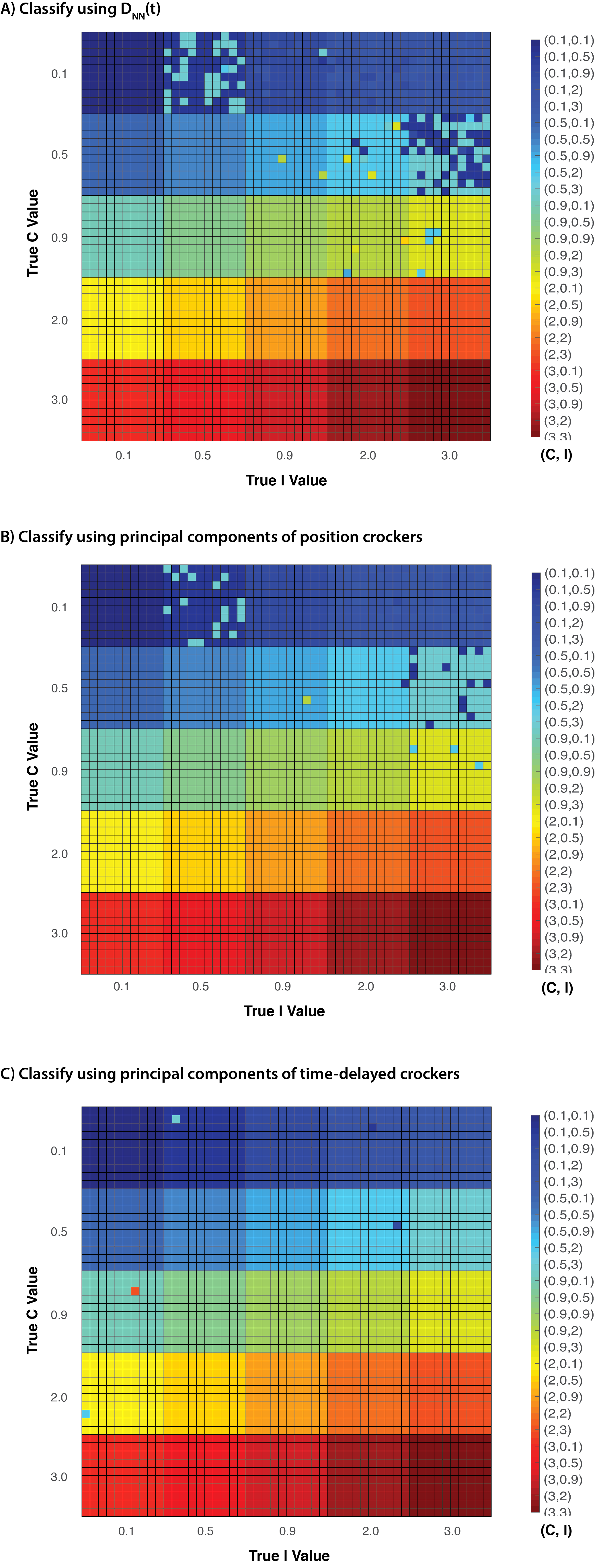}
    \caption{Results from the supervised classification of Section~\ref{sec:supervised}. In each plot, the location of a small square denotes the true underlying parameter vector ($C$,$\ell$) that generates one simulation. For instance, the $10 \times 10$ bin in the top left corner represents 100 simulations for $(C,\ell)= (0.1,0.1)$. The color of each square denotes a linear SVM's out-of-sample classification. For example, a dark blue square designates a classification of $(C,\ell)= (0.1,0.1)$. (A) Classification based on average distance to nearest neighbor, $D_{NN}(t)$. (B) Classification using the concatenation of $b_0$ and $b_1$ crockers computed from 2-D position data and reduced down to 87 dimensions via PCA for fair comparison with (A). (C) Like (B), but based on 4-D position data incorporating time delay.}
    \label{fig:supervised_plots}
\end{figure}

Fig.~\ref{fig:supervised_plots} depicts the
out-of-sample parameter classification results from linear SVMs, as described in Section~\ref{sec:supervised}. Note that in this figure, we are depicting the classification of each individual simulation and not binning these classifications together as we do in the confusion matrices of Fig.~\ref{fig:unsup_confusion}. Because of the high supervised classification accuracies, depicting the individual classifications is more informative than the summary confusion matrix.  In panel (A), the linear SVM is trained on feature vectors comprised of $D_{NN}(t)$ without any dimensionality reduction. We observe a high misclassification from parameters $(C,\ell)=(0.1,0.5)$ and $(0.5,3.0)$ as each other, as well as simulations from $(C,\ell)=(0.1,0.9),(0.1,2.0),(0.1,3.0)$. All of these parameter choices produce the collective swarm phenotype (see Table~\ref{tab:phenotype_params}), suggesting that this is the most difficult phenotype for parameter recovery. In panel (B), the linear SVM is trained on feature vectors comprised of the concatenation of $b_0$ and $b_1$ crockers derived from 2-D position data with dimensionality reduction down to 87 (to match the $D_{NN}(t)$ dimensionality). Here, we observe a marked reduction in the misclassification of the collective swarm parameter values as compared to Panel (A). In panel (C), the linear SVM is trained on feature vectors comprised of $b_0$ and $b_1$ crockers derived from 4-D time-delayed data with dimensionality reduction down to 87. Here, we observe very accurate classifications with only seven simulations being misclassified out of 2500.



\nocite{*}
\bibliography{references.bib}

\end{document}